\documentclass[11pt]{article}
\usepackage{amssymb,amsfonts,amsmath, psfrag,eepic,}
\newtheorem{theo}{Theorem}

\newtheorem{lem} [theo]{Lemma}

\makeatletter \@addtoreset{equation}{section}
\@addtoreset{theo}{section} \makeatother

\def\qed{\hfill \rule{4pt}{7pt}}
\def\pf{\noindent {\it Proof.} }

\begin{document}

\begin{center}
{\LARGE\bf Parity Reversing Involutions on Plane Trees

and 2-Motzkin Paths}\\[12pt]

William Y.C. Chen$^{*,1}$, Louis W. Shapiro$^{\dag,2}$, Laura L.M. Yang$^{*,3}$ \\[6pt]
$^*$Center for Combinatorics, LPMC\\
Nankai University,
Tianjin 300071, P. R. China\\[5pt]
$^\dag$Department of Mathematics, Howard University, Washington, DC 20059, USA\\[5pt]

$^1$chen@nankai.edu.cn, $^2$lshapiro@howard.edu,
$^3$yanglm@hotmail.com
\end{center}

\vskip 6mm

\noindent {\bf Abstract.} The problem of counting plane trees with
$n$ edges and an even or an odd number of leaves was studied by
Eu, Liu and Yeh, in connection with an identity on coloring nets
due to Stanley. This identity was also obtained by Bonin, Shapiro
and Simion in their study of Schr\"oder paths,  and it was
recently derived by Coker using the Lagrange inversion formula. An
equivalent problem for partitions was independently studied by
Klazar. We present three parity reversing involutions, one for
unlabelled plane trees, the other for labelled plane trees and one
for $2$-Motzkin paths which are in one-to-one correspondence with
Dyck paths.

 \vskip 8mm



\noindent
 {\bf Corresponding Author:} William Y. C. Chen, Email:
 chen@nankai.edu.cn

\vskip 1cm

\section{Introduction}

The set of unlabelled (rooted) plane trees with $n$ edges is
denoted by  $\mathcal{P}_n$, and is counted by the Catalan number
\[ c_n = \frac{1}{n + 1}{2n \choose n}.\]
The reader is referred to the survey of Stanley
\cite{stanley,stanleyp} and references therein  for combinatorial
objects enumerated by the Catalan numbers.  In a recent paper by
Eu, Liu and Yeh \cite{ely}, the authors consider the problem of
counting plane trees with further specification on the parity of
the number of leaves.  Let $P_e(n)$ ($P_o(n)$, respectively)
denote the number of plane trees with $n$ edges and an even (odd,
respectively) number of leaves. Eu, Liu and Yeh obtain the
following result by using generating functions.

\begin{theo}[Eu-Liu-Yeh]\label{main} The following relations hold,
\begin{eqnarray}
P_e(2n)~~~~-~~~P_o(2n)~~~ & = & 0 \label{even}\\
P_e(2n+1)~-~P_o(2n+1)&=&(-1)^{n+1}c_n. \label{odd}
\end{eqnarray}
\end{theo}

Clearly, from the above theorem one can express $P_e(n)$ and
$P_o(n)$ in terms of the Catalan numbers. Note that  this identity
was also obtained by Bonin, Shapiro and Simion \cite{bss} in their
study of Schr\"oder paths,  and it was recently derived by Coker
\cite{coker} using the Lagrange inversion formula. Klazar
\cite{klazar} also uses generating functions to derive equivalent
results for set partitions with restrictions on the parity of the
number of blocks. In fact, plane trees with $n$ edges and $k$
leaves are in one-to-one correspondence with noncrossing
partitions of $\{1,\ldots,n\}$ with $k$ blocks \cite{dz,
prodinger}.

Two combinatorial proofs of the relation (\ref{even}) are given in
\cite{ely}.  The relation (\ref{odd}) is analogous to an identity
of Stanley \cite{stanleyp} on coloring nets. It is shown in
\cite{ely} that (\ref{odd}) is equivalent to that of Stanley by a
correspondence of Deutsch \cite{deutsch} on equidistribution of
the number of even-level vertices and the number of leaves on the
set of plane trees with $n$ edges.

The objective of this paper is to give three parity reversing
involutions for both the relations (\ref{even}) and (\ref{odd}),
the first involution is based on unlabelled plane trees, and the
second is based on labelled plane trees and a decomposition
algorithm in \cite{chen}. The last involution is based on a
bijection between $2$-Motzkin paths and unlabelled plane trees in
\cite{ds2002}.

\section{An involution on plane trees}

We begin with an observation that for any plane tree with $n$
edges one may attach to each vertex a leaf as its first child to
form a plane tree with $2n+1$ edges. We use $\mathcal{B}_{n}$ to
denote the set of plane trees with $n$ edges such that any leaf is
the first child of some internal vertex and the first child of any
internal vertex must be a leaf. Notice that $\mathcal{B}_n$ is
empty if $n$ is even and
\[
|\mathcal{B}_{2n+1}|=c_n.
\]

Our involution is based on the set
$\mathcal{P}_n\setminus\mathcal{B}_n$. We define the parity of a
plane tree as the parity of the number of leaves. Moreover, we
define the sign of a plane tree as $-1$ if it is odd, and as 1 if
it is even. For any nonroot vertex $v$ of a plane tree, we say
that $v$ is \emph{legal} if $v$ is an internal vertex, but it is
not the first child of some internal vertex, or if $v$ is the
first leaf child of some internal vertex. Otherwise, $v$ is called
\emph{illegal}. A plane tree $T$ is said to be \emph{legal} if
every nonroot vertex of $T$ is legal. In particular, the plane
tree with only one vertex is illegal. In other words,
$\mathcal{B}_n$ is the set of legal trees with $n$ edges.

\begin{theo}
There is a parity reversing involution $\Phi$ on the set
$\mathcal{P}_{n}\setminus\mathcal{B}_{n}$.
\end{theo}

\pf The involution can be described recursively. Let $T$ be a
plane tree in $\mathcal{P}_{n}\setminus \mathcal{B}_{n}$. We now
conduct a depth first search for an illegal vertex of $T$ in the
following order: Let $v_1, v_2, \ldots, v_k$ be the children of
the root of $T$ from left to right and $T_i$ be the subtree of $T$
rooted at the vertex $v_i$ for $1\leq i\leq k$. Then we search for
an illegal vertex in $T_k$. If $T_k$ is legal, then we conduct the
search for $T_{k-1}$, and so on. If $T_2, \ldots, T_k$ are all
legal, then the first child $v_1$ of $T$ must be an internal
vertex which implies that $v_1$ is illegal. Using the above search
scheme, we can find an illegal vertex $v$ of $T$ which is the
first vertex encountered while implementing the above search.

Let $u$ be the father of $v$ and $T_u$ be the subtree of $T$
rooted at the vertex $u$. We now have two cases. (1) The vertex
$v$ is a leaf, but it is not the first child of $u$. (2) The
vertex $v$ is an internal vertex, and it is the first child of
$u$. In this case, all the subtrees rooted at the other children
of $u$ are legal.

For Case (1), let $w_1, w_2, \ldots, w_i$ be the children of $u$
that are to the left of $v$. We now cut off the edges between $u$
and $w_1, w_2, \ldots, w_i$, and move the subtrees $T_{w_1} ,
\ldots, T_{w_i}$ as subtrees of $v$ in the same order. Let
$\Phi(T)$ denote the resulted tree.
 Note that in the search process for $\Phi(T)$,
 the vertex $v$ is still the first encountered illegal vertex.

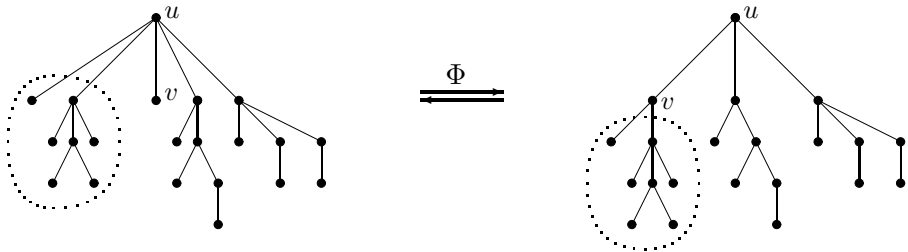
\begin{figure}[h,t]
\begin{center}
\begin{picture}(340,100)
\setlength{\unitlength}{1.1mm} \linethickness{0.4pt}

\put(25,15){\circle*{1}} \put(15,25){\circle*{1}}
\put(15,25){\line(1,-1){10}} \put(25,15){\line(2,-1){10}}
\put(35,10){\circle*{1}} \put(25,15){\line(1,-1){5}}
\put(30,10){\circle*{1}} \put(25,15){\line(0,-1){5}}
\put(25,10){\circle*{1}} \put(35,10){\line(0,-1){5}}
\put(35,5){\circle*{1}} \put(30,10){\line(0,-1){5}}
\put(30,5){\circle*{1}}

\put(15,25){\line(0,-1){10}} \put(15,15){\circle*{1}}

\put(20,15){\line(0,-1){5}} \put(20,15){\line(-1,-2){2.5}}
\put(20,15){\circle*{1}} \put(20,10){\circle*{1}}
\put(17.5,10){\circle*{1}}

\put(20,10){\line(-1,-2){2.5}} \put(20,10){\line(1,-2){2.5}}
\put(22.5,5){\line(0,-1){5}} \put(17.5,5){\circle*{1}}
\put(22.5,5){\circle*{1}} \put(22.5,0){\circle*{1}}

\put(15,25){\line(1,-2){5}} \put(20,15){\circle*{1}}

\put(15,25){\line(-1,-1){10}} \put(5,15){\line(1,-2){2.5}}
\put(5,15){\line(-1,-2){2.5}} \put(5,15){\line(0,-1){5}}
\put(5,15){\circle*{1}} \put(7.5,10){\circle*{1}}
\put(2.5,10){\circle*{1}} \put(5,10){\circle*{1}}
\put(5,10){\line(-1,-2){2.5}} \put(5,10){\line(1,-2){2.5}}
\put(2.5,5){\circle*{1}} \put(7.5,5){\circle*{1}}

\put(15,25){\line(-3,-2){15}} \put(0,15){\circle*{1}}

\linethickness{0.7pt}

\qbezier[10](2.5,2)(-3,2)(-3,10)
\qbezier[10](2.5,2)(10.5,2)(10.5,10)
\qbezier[10](-3,10)(-3,18)(2.5,18)
\qbezier[10](2.5,18)(10.5,18)(10.5,10)
\linethickness{1pt} \put(47,16){\vector(1,0){10}}
\put(57,15){\vector(-1,0){10}} \put(50,17){$\Phi$}
\linethickness{0.4pt} \put(75,15){\line(1,1){10}}
\put(85,25){\circle*{1}} \put(75,15){\line(-1,-1){5}}
\put(75,15){\circle*{1}} \put(70,10){\circle*{1}}
\put(75,10){\line(0,1){5}} \put(75,10){\line(1,-2){2.5}}
\put(75,10){\line(-1,-2){2.5}} \put(75,10){\line(0,-1){5}}
\put(75,10){\circle*{1}} \put(77.5,5){\circle*{1}}
\put(72.5,0){\circle*{1}} \put(75,5){\circle*{1}}
\put(75,5){\line(-1,-2){2.5}} \put(75,5){\line(1,-2){2.5}}
\put(72.5,5){\circle*{1}} \put(77.5,0){\circle*{1}}
\put(85,25){\line(0,-1){10}} \put(85,15){\line(1,-2){2.5}}
\put(85,15){\line(-1,-2){2.5}} \put(85,15){\circle*{1}}
\put(87.5,10){\circle*{1}} \put(82.5,10){\circle*{1}}

\put(87.5,10){\line(-1,-2){2.5}} \put(87.5,10){\line(1,-2){2.5}}
\put(90,5){\line(0,-1){5}} \put(85,5){\circle*{1}}
\put(90,5){\circle*{1}} \put(90,0){\circle*{1}}
\put(85,25){\circle*{1}} \put(85,25){\line(1,-1){10}}
\put(95,15){\circle*{1}} \put(95,15){\line(2,-1){10}}
\put(105,10){\circle*{1}} \put(95,15){\line(1,-1){5}}
\put(100,10){\circle*{1}} \put(95,15){\line(0,-1){5}}
\put(95,10){\circle*{1}} \put(105,10){\line(0,-1){5}}
\put(105,5){\circle*{1}} \put(100,10){\line(0,-1){5}}
\put(100,5){\circle*{1}}
\linethickness{0.7pt} \qbezier[10](74,-3)(67,-3)(67,5)
\qbezier[10](67,5)(67,13)(74,13)
\qbezier[10](74,13)(80.5,13)(80.5,5)
\qbezier[10](80.5,5)(80.5,-3)(74,-3)

\linethickness{0.4pt} \put(16,25){\small{$u$}}
\put(16,15){\small{$v$}} \put(86,25){\small{$u$}}
\put(76,14){\small{$v$}}
\end{picture}
\end{center}
\caption{Involution $\Phi$ on plane trees.} \label{phi}
\end{figure}

For Case (2), we may reverse the construction for Case (1). Hence
we obtain a parity reversing involution $\Phi$. See Fig.
\ref{phi}.\qed

Note that $\mathcal{B}_{2n}$ is empty and any plane tree in
$\mathcal{B}_{2n+1}$ has $n+1$ leaves.  The involution $\Phi$
implies a combinatorial proof of relations (\ref{even}) and
(\ref{odd}) if the signs of the plane trees are taken into
account. We also note that the involution $\Phi$ is different from
the involution of Eu-Liu-Yeh \cite{ely} for the relation
(\ref{even}).

\section{A bijective algorithm for labelled trees}

In this section, we give a parity reversing involution on labelled
plane trees that also leads to a combinatorial interpretation of
the relations for labelled plane trees. Note that the number of
labelled plane tress with $n$ edges, or $n+1$ vertices,  equals
$(n+1)!$ times the number of unlabelled plane trees with $n$
edges, which we denote by
\[   t_n = (n+1)!\, c_n = {(2n)!\over n! } .\]

In  \cite{chen}, the author gives a  bijective algorithm to
decompose a labelled plane tree on $\{1, 2, \ldots, n+1\}$ into a
set $F$ of $n$ matches with labels $\{1,\ldots,n,n+1, (n+2)^*,
\ldots, (2n)^*\}$, where a match is a rooted tree with two
vertices. The reverse procedure of the decomposition algorithm is
the following merging algorithm. We start with a set $F$ of
matches on $\{1, \ldots, n+1, (n+2)^*, \ldots, (2n)^*\}$. A vertex
labelled by a mark $*$ is called a marked vertex.

\begin{figure}[h,t]
\begin{center}
\begin{picture}(420,50)
\setlength{\unitlength}{1.1mm} \linethickness{0.4pt}

\put(5,10){\line(-1,-1){5}} \put(0,5){\circle*{1}}
\put(0,5){\line(0,-1){5}} \put(0,0){\circle*{1}}
\put(5,10){\circle*{1}} \put(5,5){\circle*{1}}
\put(5,0){\line(0,1){10}} \put(5,0){\circle*{1}} \put(3,13){$i$}

\put(11,4){\large\bf$+$} \put(20,0){\circle*{1}}
\put(20,10){\circle*{1}} \put(20,0){\line(0,1){10}}
\put(18,13){$j^*$} \put(25,4){\large\bf$\Rightarrow$}

\put(35,0){\circle*{1}} \put(35,5){\circle*{1}}
\put(35,0){\line(0,1){5}} \put(40,10){\circle*{1}}
\put(35,5){\line(1,1){5}} \put(40,5){\circle*{1}}
\put(40,10){\line(0,-1){10}} \put(40,0){\circle*{1}}
\put(40,10){\line(1,-1){5}} \put(45,5){\circle*{1}}

\put(39,13){$i$} \put(2,-7){$T$} \put(18,-7){$T^*$}

\put(80,0){\circle*{1}} \put(80,0){\line(1,2){5}}
\put(85,10){\circle*{1}} \put(90,0){\line(-1,2){5}}
\put(90,0){\circle*{1}} \put(83,13){$i$}

\put(95,4){\large\bf$+$}

\put(105,0){\circle*{1}} \put(105,10){\circle*{1}}
\put(105,0){\line(0,1){10}} \put(106,-1){$j^*$}
\put(110,4){\large\bf$\Rightarrow$} \put(120,0){\circle*{1}}
\put(125,5){\circle*{1}} \put(120,0){\line(1,1){5}}
\put(130,0){\circle*{1}} \put(130,0){\line(-1,1){5}}
\put(125,10){\circle*{1}} \put(125,5){\line(0,1){5}}
\put(127,4){$i$} \put(85,-7){$T$} \put(103,-7){$T^*$}
\end{picture}
\end{center}
\caption{A horizontal merge and vertical merge.}\label{merge}
\end{figure}
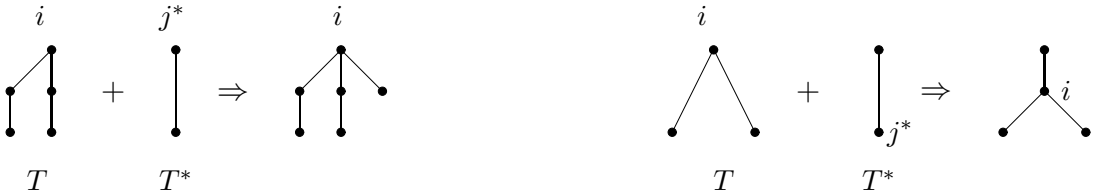

\noindent (1) Find the tree $T$ with the smallest root in which no
vertex is marked. Let $i$ be the root of
$T$.\\
(2) Find the tree $T^*$ in $F$ that contains the smallest marked
vertex. Let $j^*$ be this marked vertex.\\
(3) If $j^*$ is the root of $T^*$, then merge $T$ and $T^*$ by
identifying $i$ and $j^*$, keep $i$ as the new vertex, and put the
subtrees of $T^*$ to the right of $T$. The operation is called a
\emph{horizontal merge}. If $j^*$ is a leaf of $T^*$, then replace
$j^*$ with $T$ in $T^*$. This operation is called a \emph{vertical
merge}. See Figure \ref{merge}.

\noindent (4) Repeat the above procedure until $F$ becomes a
labelled tree.

For any set $F$ of $n$ matches labelled by $\{1, \ldots, n+1,
(n+2)^*, \ldots, (2n)^*\}$, a match is said to be \emph{pure} if
 it consists of either two unmarked vertices or two marked vertices.
 We use $\mathcal{A}_n$ to denote the set of
pure matches on $\{1, \ldots, n+1, (n+2)^*, \ldots, (2n)^*\}$. It
is easy to see that $|\mathcal{A}_{2n}|=0$ and
\[ |\mathcal{A}_{2n+1} | =   {(2n)! \over n!} \, { (2n+2)! \over (n+1)!}  = t_n t_{n+1}\, .\]

We are now ready to give an involution for the following
 labelled version of Theorem \ref{main}. Let $Q_e(n)$ ($Q_o(n)$,
respectively) be the number of labelled trees with $n$ edges and
an even (odd, respectively) number of leaves.
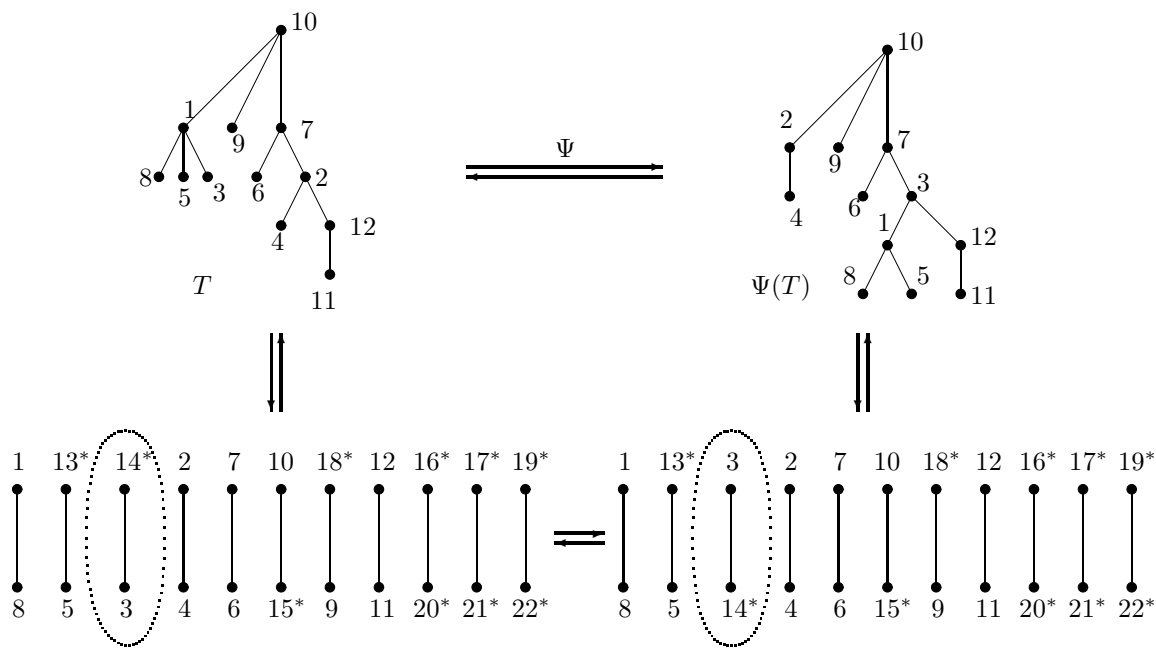
\begin{figure}[h,t]
\begin{center}
\begin{picture}(400,220)
\setlength{\unitlength}{1.3mm} \linethickness{0.4pt}
\put(-2,0){\line(0,1){10}} \put(-2,0){\circle*{1}}
\put(-2,10){\circle*{1}} \put(-2.5,12){\small{1}}
\put(-2.5,-3){\small{8}}

\put(3,10){\line(0,-1){10}} \put(3,10){\circle*{1}}
\put(3,0){\circle*{1}} \put(1.6,12){\small{$13^*$}}
\put(2.5,-3){\small{5}}

\put(9,10){\line(0,-1){10}} \put(9,10){\circle*{1}}
\put(9,0){\circle*{1}} \put(8,12){\small{$14^*$}}
\put(8.5,-3){\small{3}}

\linethickness{0.7pt} \qbezier[20](9,-6)(5,-5)(5,5)
\qbezier[20](9,-6)(13,-6)(13,5) \qbezier[20](5,5)(5,16)(9,16)
\qbezier[20](9,16)(13,16)(13,5)

\linethickness{0.4pt} \put(15,10){\line(0,-1){10}}
\put(15,10){\circle*{1}} \put(15,0){\circle*{1}}
\put(14.5,12){\small{$2$}} \put(14.5,-3){\small{4}}

\put(20,10){\line(0,-1){10}} \put(20,10){\circle*{1}}
\put(20,0){\circle*{1}} \put(19.5,12){\small{$7$}}
\put(19.5,-3){\small{6}}

\put(25,10){\line(0,-1){10}} \put(25,10){\circle*{1}}
\put(25,0){\circle*{1}} \put(23.6,12){\small{$10$}}
\put(23.6,-3){\small{$15^*$}}

\put(30,10){\line(0,-1){10}} \put(30,10){\circle*{1}}
\put(30,0){\circle*{1}} \put(28.6,12){\small{$18^*$}}
\put(29.5,-3){\small{$9$}}

\put(35,10){\line(0,-1){10}} \put(35,10){\circle*{1}}
\put(35,0){\circle*{1}} \put(34,12){\small{$12$}}
\put(34,-3){\small{$11$}}

\put(40,10){\line(0,-1){10}} \put(40,10){\circle*{1}}
\put(40,0){\circle*{1}} \put(38.5,12){\small{$16^*$}}
\put(38.5,-3){\small{$20^*$}}

\put(45,10){\line(0,-1){10}} \put(45,10){\circle*{1}}
\put(45,0){\circle*{1}} \put(43.6,12){\small{$17^*$}}
\put(43.5,-3){\small{$21^*$}}

\put(50,10){\line(0,-1){10}} \put(50,10){\circle*{1}}
\put(50,0){\circle*{1}} \put(48.6,12){\small{$19^*$}}
\put(48.6,-3){\small{$22^*$}}
\linethickness{1pt} \put(53,5.5){\vector(1,0){5}}
\put(58,4.5){\vector(-1,0){5}}
\linethickness{0.4pt} \put(60,0){\line(0,1){10}}
\put(60,0){\circle*{1}} \put(60,10){\circle*{1}}
\put(59.5,12){\small{1}} \put(59.5,-3){\small{8}}

\put(65,10){\line(0,-1){10}} \put(65,10){\circle*{1}}
\put(65,0){\circle*{1}} \put(63.6,12){\small{$13^*$}}
\put(64.5,-3){\small{5}}

\put(71,10){\line(0,-1){10}} \put(71,10){\circle*{1}}
\put(71,0){\circle*{1}} \put(70,-3){\small{$14^*$}}
\put(70.5,12){\small{3}}

\linethickness{0.7pt} \qbezier[20](71,-6)(67,-6)(67,5)
\qbezier[20](71,-6)(75,-6)(75,5) \qbezier[20](67,5)(67,16)(71,16)
\qbezier[20](71,16)(75,16)(75,5) \linethickness{0.4pt}

\put(77,10){\line(0,-1){10}} \put(77,10){\circle*{1}}
\put(77,0){\circle*{1}} \put(76.5,12){\small{$2$}}
\put(76.5,-3){\small{4}}

\put(82,10){\line(0,-1){10}} \put(82,10){\circle*{1}}
\put(82,0){\circle*{1}} \put(81.5,12){\small{$7$}}
\put(81.5,-3){\small{6}}

\put(87,10){\line(0,-1){10}} \put(87,10){\circle*{1}}
\put(87,0){\circle*{1}} \put(85.6,12){\small{$10$}}
\put(85.6,-3){\small{$15^*$}}

\put(92,10){\line(0,-1){10}} \put(92,10){\circle*{1}}
\put(92,0){\circle*{1}} \put(90.6,12){\small{$18^*$}}
\put(91.5,-3){\small{$9$}}

\put(97,10){\line(0,-1){10}} \put(97,10){\circle*{1}}
\put(97,0){\circle*{1}} \put(96,12){\small{$12$}}
\put(96,-3){\small{$11$}}

\put(102,10){\line(0,-1){10}} \put(102,10){\circle*{1}}
\put(102,0){\circle*{1}} \put(100.5,12){\small{$16^*$}}
\put(100.5,-3){\small{$20^*$}}

\put(107,10){\line(0,-1){10}} \put(107,10){\circle*{1}}
\put(107,0){\circle*{1}} \put(105.6,12){\small{$17^*$}}
\put(105.5,-3){\small{$21^*$}}

\put(112,10){\line(0,-1){10}} \put(112,10){\circle*{1}}
\put(112,0){\circle*{1}} \put(110.6,12){\small{$19^*$}}
\put(110.6,-3){\small{$22^*$}}
\linethickness{1pt} \put(25,18){\vector(0,1){8}}
\put(24,26){\vector(0,-1){8}} \put(85,18){\vector(0,1){8}}
\put(84,26){\vector(0,-1){8}}
\linethickness{0.4pt} \put(16,30){\small$T$}

\put(25,57){\line(-1,-2){5}} \put(20,47){\circle*{1}}
\put(20,44.5){\small{9}}

\put(25,47){\line(1,-2){2.5}} \put(25,47){\line(-1,-2){2.5}}
\put(25,47){\circle*{1}} \put(27.5,42){\circle*{1}}
\put(22.5,42){\circle*{1}} \put(27,46){\small{7}}
\put(28.5,41){\small{2}} \put(22,39.5){\small{6}}

\put(27.5,42){\line(-1,-2){2.5}} \put(27.5,42){\line(1,-2){2.5}}
\put(30,37){\line(0,-1){5}} \put(25,37){\circle*{1}}
\put(24,34.5){\small{4}} \put(30,37){\circle*{1}}
\put(32,36){\small{12}} \put(30,32){\circle*{1}}
\put(28,28.5){\small{11}}

\put(25,57){\line(0,-1){10}} \put(25,57){\circle*{1}}
\put(26,57){\small{10}}

\put(25,57){\line(-1,-1){10}} \put(15,47){\line(1,-2){2.5}}
\put(15,47){\line(-1,-2){2.5}} \put(15,47){\line(0,-1){5}}
\put(15,47){\circle*{1}} \put(15,48){\small{1}}
\put(17.5,42){\circle*{1}}
\put(14.5,39){\small{5}}\put(12.5,42){\circle*{1}}
\put(10.5,41){\small{8}} \put(15,42){\circle*{1}}
\put(18,39.5){\small{3}}
\linethickness{1pt} \put(44,43){\vector(1,0){20}}
\put(64,42){\vector(-1,0){20}} \put(53,44){\small{$\Psi$}}
\linethickness{0.4pt}

\put(73,30){\small$\Psi(T)$} \put(77,40){\line(0,1){5}}
\put(77,40){\circle*{1}} \put(77,37){\small{4}}
\put(77,45){\circle*{1}} \put(76,47){\small{2}}

\put(77,45){\line(1,1){10}} \put(87,55){\circle*{1}}
\put(88,55){\small{10}}

\put(87,55){\line(-1,-2){5}} \put(82,45){\circle*{1}}
\put(81,42.5){\small{9}}

\put(87,55){\line(0,-1){10}} \put(87,45){\circle*{1}}
\put(88,45){\small{7}}

\put(87,45){\line(-1,-2){2.5}} \put(84.5,40){\circle*{1}}
\put(83,38){\small{6}}

\put(87,45){\line(1,-2){2.5}} \put(89.5,40){\circle*{1}}
\put(90,40.5){\small{3}}

\put(89.5,40){\line(-1,-2){2.5}} \put(87,35){\circle*{1}}
\put(86,36.5){\small{1}}

\put(87,35){\line(-1,-2){2.5}} \put(84.5,30){\circle*{1}}
\put(82.5,31){\small{8}}

\put(87,35){\line(1,-2){2.5}} \put(89.5,30){\circle*{1}}
\put(90,31){\small{5}}

\put(89.5,40){\line(1,-1){5}} \put(94.5,35){\circle*{1}}
\put(95.5,35){\small{12}}

\put(94.5,35){\line(0,-1){5}} \put(94.5,30){\circle*{1}}
\put(95.5,29){\small{11}}
\end{picture}
\end{center}
\caption{Involution $\Psi$ on labelled plane trees.} \label{psi}
\end{figure}

\begin{theo}
The following relations hold,
\begin{eqnarray*}
Q_e(2n)~~~~~ -  ~~Q_o(2n)~~~ & = & 0\\
Q_e(2n+1) ~~-~~ Q_o(2n+1)&=&(-1)^{n+1} |\mathcal{A}_{2n+1}|.
\end{eqnarray*}
\end{theo}

\pf We define the sign of a labelled plane tree in the same way as
in the unlabelled case. The involution $\Psi$ is built on the set
of labelled plane trees whose match decompositions contain a match
with mixed vertices (one marked, the other unmarked). Given such a
plane tree $T$, we decompose it into matches. Then we choose the
match with mixed vertices such that the unmarked vertex is
minimum. The involution is simply to turn the chosen match up side
down, see Figure \ref{psi}. Note that the number of leaves of a
plane tree $T$ equals the number of unmarked leaves of the matches
in the corresponding decomposition, see \cite{chen}. \qed

\section{An involution on 2-Motzkin paths}

The notion of $2$-Motzkin paths has proved to be very useful
representation of several combinatorial objects such as Dyck
paths, plane trees, noncrossing partitions \cite{ds2002}. For the
purpose of this paper, we need the bijection between $2$-Motzkin
paths and plane trees. Recall that a {\it $2$-Motzkin path} is a
lattice path starting and ending on the horizontal axis but never
going below it, with up steps $(1,1)$, level steps $(1,0)$, and
down steps $(1,-1)$, where the level steps can be either of two
kinds: straight and wavy. The {\it length} of a path is defined to
be the number of its steps. The set of $2$-Motzkin paths of length
$n-1$ is denoted by $\mathcal{M}_n$. The following bijection is
given in \cite{ds2002}.

\begin{lem} \label{motzkin-tree}
There is a bijection between plane trees with $n$ edges and
$2$-Motzkin paths of length $n-1$, such that the number of leaves
of a plane tree minus $1$ equals the sum of the number of up steps
and the number of wavy level steps in the corresponding
$2$-Motzkin path.
\end{lem}

Note that the set of $2$-Motzkin paths without level steps reduces
to the set of Dyck paths \cite{stanley}. We use $\mathcal{D}_n$ to
denote the set of $2$-Motzkin paths of length $n-1$ without level
steps, namely, the set of Dyck paths of length $n-1$. Obviously,
$\mathcal{D}_n$ is empty if $n$ is even and
$$
|\mathcal{D}_{2n+1}|=c_n.
$$

We obtain the following involution on the set of $2$-Motzkin paths
which gives the third combinatorial interpretation of the
relations (\ref{even}) and (\ref{odd}). Note that the parity of a
$2$-Motzkin path is meant to be the parity  of one plus the sum of
the number of up steps and the number of wavy steps, as indicated
by the above Lemma \ref{motzkin-tree}.

\begin{theo}
There is a parity reversing involution $\Upsilon$ on
$\mathcal{M}_n\setminus\mathcal{D}_n$.
\end{theo}

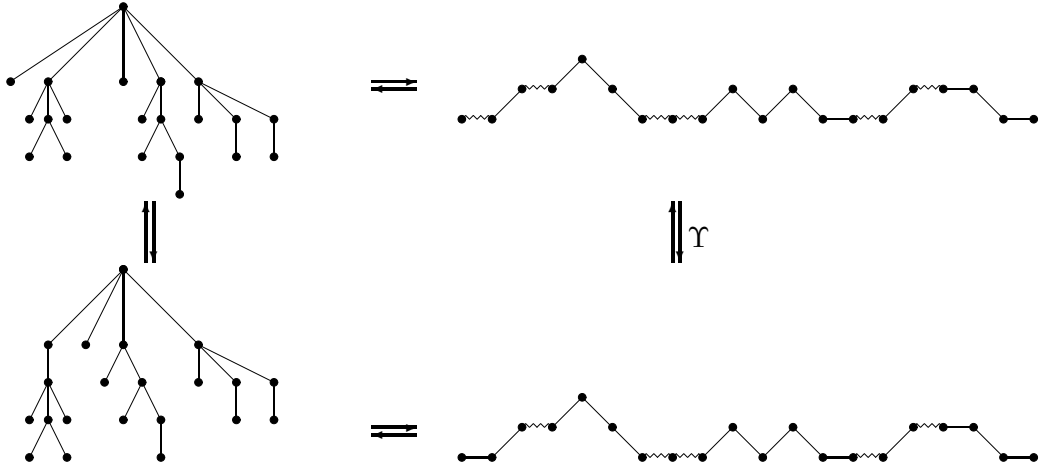
\begin{figure}[h,t]
\begin{center}
\begin{picture}(400,200)
\setlength{\unitlength}{1mm} \linethickness{0.4pt}

\put(25,50){\circle*{1}} \put(15,60){\circle*{1}}
\put(15,60){\line(1,-1){10}} \put(25,50){\line(2,-1){10}}
\put(35,45){\circle*{1}} \put(25,50){\line(1,-1){5}}
\put(30,45){\circle*{1}} \put(25,50){\line(0,-1){5}}
\put(25,45){\circle*{1}} \put(35,45){\line(0,-1){5}}
\put(35,40){\circle*{1}} \put(30,45){\line(0,-1){5}}
\put(30,40){\circle*{1}}

\put(15,60){\line(0,-1){10}} \put(15,50){\circle*{1}}

\put(20,50){\line(0,-1){5}} \put(20,50){\line(-1,-2){2.5}}
\put(20,50){\circle*{1}} \put(20,45){\circle*{1}}
\put(17.5,45){\circle*{1}}

\put(20,45){\line(-1,-2){2.5}} \put(20,45){\line(1,-2){2.5}}
\put(22.5,40){\line(0,-1){5}} \put(17.5,40){\circle*{1}}
\put(22.5,40){\circle*{1}} \put(22.5,35){\circle*{1}}

\put(15,60){\line(1,-2){5}} \put(20,50){\circle*{1}}

\put(15,60){\line(-1,-1){10}} \put(5,50){\line(1,-2){2.5}}
\put(5,50){\line(-1,-2){2.5}} \put(5,50){\line(0,-1){5}}
\put(5,50){\circle*{1}} \put(7.5,45){\circle*{1}}
\put(2.5,45){\circle*{1}} \put(5,45){\circle*{1}}
\put(5,45){\line(-1,-2){2.5}} \put(5,45){\line(1,-2){2.5}}
\put(2.5,40){\circle*{1}} \put(7.5,40){\circle*{1}}

\put(15,60){\line(-3,-2){15}} \put(0,50){\circle*{1}}
\linethickness{1pt} \put(18,26){\vector(0,1){8}}
\put(19,34){\vector(0,-1){8}}
\linethickness{0.4pt} \put(5,15){\line(1,1){10}}
\put(15,25){\circle*{1}} \put(15,25){\line(-1,-2){5}}
\put(5,15){\circle*{1}} \put(10,15){\circle*{1}}
\put(5,10){\line(0,1){5}} \put(5,10){\line(1,-2){2.5}}
\put(5,10){\line(-1,-2){2.5}} \put(5,10){\line(0,-1){5}}
\put(5,10){\circle*{1}} \put(7.5,5){\circle*{1}}
\put(2.5,0){\circle*{1}} \put(5,5){\circle*{1}}
\put(5,5){\line(-1,-2){2.5}} \put(5,5){\line(1,-2){2.5}}
\put(2.5,5){\circle*{1}} \put(7.5,0){\circle*{1}}
\put(15,25){\line(0,-1){10}} \put(15,15){\line(1,-2){2.5}}
\put(15,15){\line(-1,-2){2.5}} \put(15,15){\circle*{1}}
\put(17.5,10){\circle*{1}} \put(12.5,10){\circle*{1}}

\put(17.5,10){\line(-1,-2){2.5}} \put(17.5,10){\line(1,-2){2.5}}
\put(20,5){\line(0,-1){5}} \put(15,5){\circle*{1}}
\put(20,5){\circle*{1}} \put(20,0){\circle*{1}}
\put(15,25){\circle*{1}} \put(15,25){\line(1,-1){10}}
\put(25,15){\circle*{1}} \put(25,15){\line(2,-1){10}}
\put(35,10){\circle*{1}} \put(25,15){\line(1,-1){5}}
\put(30,10){\circle*{1}} \put(25,15){\line(0,-1){5}}
\put(25,10){\circle*{1}} \put(35,10){\line(0,-1){5}}
\put(35,5){\circle*{1}} \put(30,10){\line(0,-1){5}}
\put(30,5){\circle*{1}}
\linethickness{1pt} \put(48,4){\vector(1,0){6}}
\put(54,3){\vector(-1,0){6}}
\linethickness{0.4pt} \put(60,0){\circle*{1}}
\put(60,0){\line(1,0){4}}

\put(64,0){\circle*{1}} \put(64,0){\line(1,1){4}}

\put(68,4){\circle*{1}} \multiput(68,4)(1,0){4}{\line(1,1){0.5}}
\multiput(68.5,4.5)(1,0){4}{\line(1,-1){0.5}}

\put(72,4){\circle*{1}} \put(72,4){\line(1,1){4}}

\put(76,8){\circle*{1}} \put(76,8){\line(1,-1){8}}
\put(80,4){\circle*{1}}

\put(84,0){\circle*{1}} \multiput(84,0)(1,0){8}{\line(1,1){0.5}}
\multiput(84.5,0.5)(1,0){8}{\line(1,-1){0.5}}
\put(88,0){\circle*{1}}

\put(92,0){\circle*{1}} \put(92,0){\line(1,1){4}}
\put(96,4){\circle*{1}} \put(96,4){\line(1,-1){4}}

\put(100,0){\circle*{1}} \put(100,0){\line(1,1){4}}
\put(104,4){\circle*{1}} \put(104,4){\line(1,-1){4}}

\put(108,0){\circle*{1}} \put(108,0){\line(1,0){4}}

\put(112,0){\circle*{1}} \multiput(112,0)(1,0){4}{\line(1,1){0.5}}
\multiput(112.5,0.5)(1,0){4}{\line(1,-1){0.5}}

\put(116,0){\circle*{1}} \put(116,0){\line(1,1){4}}

\put(120,4){\circle*{1}} \multiput(120,4)(1,0){4}{\line(1,1){0.5}}
\multiput(120.5,4.5)(1,0){4}{\line(1,-1){0.5}}

\put(124,4){\circle*{1}} \put(124,4){\line(1,0){4}}

\put(128,4){\circle*{1}} \put(128,4){\line(1,-1){4}}

\put(132,0){\circle*{1}} \put(132,0){\line(1,0){4}}
\put(136,0){\circle*{1}}
\linethickness{1pt} \put(88,26){\vector(0,1){8}}
\put(89,34){\vector(0,-1){8}} \put(90,28){$\Upsilon$}
\linethickness{0.4pt} \put(60,45){\circle*{1}}
\multiput(60,45)(1,0){4}{\line(1,1){0.5}}
\multiput(60.5,45.5)(1,0){4}{\line(1,-1){0.5}}

\put(64,45){\circle*{1}} \put(64,45){\line(1,1){4}}

\put(68,49){\circle*{1}} \multiput(68,49)(1,0){4}{\line(1,1){0.5}}
\multiput(68.5,49.5)(1,0){4}{\line(1,-1){0.5}}

\put(72,49){\circle*{1}} \put(72,49){\line(1,1){4}}

\put(76,53){\circle*{1}} \put(76,53){\line(1,-1){8}}
\put(80,49){\circle*{1}}

\put(84,45){\circle*{1}} \multiput(84,45)(1,0){8}{\line(1,1){0.5}}
\multiput(84.5,45.5)(1,0){8}{\line(1,-1){0.5}}
\put(88,45){\circle*{1}}

\put(92,45){\circle*{1}} \put(92,45){\line(1,1){4}}
\put(96,49){\circle*{1}} \put(96,49){\line(1,-1){4}}

\put(100,45){\circle*{1}} \put(100,45){\line(1,1){4}}
\put(104,49){\circle*{1}} \put(104,49){\line(1,-1){4}}

\put(108,45){\circle*{1}} \put(108,45){\line(1,0){4}}

\put(112,45){\circle*{1}}
\multiput(112,45)(1,0){4}{\line(1,1){0.5}}
\multiput(112.5,45.5)(1,0){4}{\line(1,-1){0.5}}

\put(116,45){\circle*{1}} \put(116,45){\line(1,1){4}}

\put(120,49){\circle*{1}}
\multiput(120,49)(1,0){4}{\line(1,1){0.5}}
\multiput(120.5,49.5)(1,0){4}{\line(1,-1){0.5}}

\put(124,49){\circle*{1}} \put(124,49){\line(1,0){4}}

\put(128,49){\circle*{1}} \put(128,49){\line(1,-1){4}}

\put(132,45){\circle*{1}} \put(132,45){\line(1,0){4}}
\put(136,45){\circle*{1}}
\linethickness{1pt} \put(48,50){\vector(1,0){6}}
\put(54,49){\vector(-1,0){6}}

\end{picture}
\end{center}
\caption{Involution $\Upsilon$ on $2$-Motzkin paths.}
\label{Upsilon}
\end{figure}

\pf For any $2$-Motzkin path in
$\mathcal{M}_n\setminus\mathcal{D}_n$, we find the first level
step and toggle this step between wavy and straight. See Fig.
\ref{Upsilon}.\qed

It is clear that the above involution $\Upsilon$ changes the
parity of the number of  wavy  steps and keeps the number of up
steps invariant. We note that the above involution is different
from the first involution as given in Section 2.

\noindent{\bf Remark.} Chen, Deutsch and Elizalde \cite{cde}
recently found a bijection between plane trees with $n$ edges and
2-Motzkin paths of length $n-1$ such that the non-rightmost leaves
are corresponded to wavy steps. Recall that a leaf is said to be
rightmost if it is the rightmost child of its parent. Our third
involution on $2$-Motzkin paths can be used to give a
combinatorial proof for the following identities of Sun \cite{sun}
which are derived by using generating functions:
\begin{theo}[\cite{sun}]
Let $T_{n,k}$ be the  number of Dyck paths of length $2n$ with $k$
udu's. Then we have
\begin{eqnarray*}
\sum_{k~even}T_{2n,k}&=&\sum_{k~odd}T_{2n,k}, \\[2pt]
\sum_{k~even}T_{2n-1,k}&=&\sum_{k~odd}T_{2n-1,k}+c_{n-1}.
\end{eqnarray*}
\end{theo}

 \vskip 5mm

\noindent{\bf Acknowledgments.} We would like to thank the
referees for helpful comments. This work was done under the
auspices of the ¡°973¡± Project on Mathematical Mechanization, the
National Science Foundation, and the Ministry of Science and
Technology of China.

\small

\end{document}